\documentclass[conference]{IEEEtran}
\IEEEoverridecommandlockouts
\usepackage{graphicx}
\usepackage{array}
\usepackage{textcomp}
\usepackage{cite}
\usepackage{url}
\usepackage{amsmath}
\usepackage{multirow}
\usepackage{bm}
\usepackage{verbatim}

\newcommand{\splitatcommas}[1]{%
	\begingroup
	\begingroup\lccode`~=`, \lowercase{\endgroup
		\edef~{\mathchar\the\mathcode`, \penalty0 \noexpand\hspace{0pt plus 0em}}%
	}\mathcode`,="8000 #1%
	\endgroup
}

\ifCLASSINFOpdf
\else
\fi
\begin{document}
\bstctlcite{IEEEexample:BSTcontrol}	
\title{Bilevel Optimization Based Transmission Expansion Planning Considering Phase Shifting Transformer}

\author{\IEEEauthorblockN{Xiaohu Zhang\IEEEauthorrefmark{1},
		Di Shi\IEEEauthorrefmark{1},
		Zhiwei Wang\IEEEauthorrefmark{1}, 
		Zhe Yu\IEEEauthorrefmark{1},
		Xinan Wang\IEEEauthorrefmark{1},
		Desong Bian\IEEEauthorrefmark{1} and
		Kevin Tomsovic\IEEEauthorrefmark{2}}
	\IEEEauthorblockA{\IEEEauthorrefmark{1}GEIRI North America}
	\IEEEauthorblockA{\IEEEauthorrefmark{2}The Univerisity of Tennessee, Knoxville\\
		Email: xiaohu.zhang@geirina.net}
	\thanks{This project is funded by State Grid Corporation of China (SGCC) under project \textsl{Research on Spatial-Temporal Multidimensional Coordination of Energy Resources}.}}

%


\maketitle

\begin{abstract}
 In this paper, the phase shifting transformer (PST) is introduced in the Transmission Expansion Planning (TEP) problem considering significant wind power integration. The proposed planning model is formulated as a bilevel program which seeks to determine the optimal strategy for the network reinforcements and the PST locations. The objective of the upper level problem is to minimize the total consumer payment, the investment cost on transmission line and PST. The lower level problems designate the electricity market clearing conditions under different load-wind scenarios. The bilevel model is transformed into a single level mixed integer linear program (MILP) by replacing each lower level problem with its primal-dual formulation. The numerical case studies based on IEEE 24-bus system demonstrate the characteristics of the proposed model. Moreover, the simulation results show that the installation of PST adds flexibility to the TEP and facilitates the integration of wind power.

\end{abstract}

\begin{IEEEkeywords}
	Bilevel optimization, transmission expansion planning, phase shifting transformer, wind power integration.
\end{IEEEkeywords}

%
\IEEEpeerreviewmaketitle

\section*{Nomenclature}
\subsection*{Indices}
\addcontentsline{toc}{section}{Nomenclature}
\begin{IEEEdescription}[\IEEEusemathlabelsep\IEEEsetlabelwidth{$V_1,V_2$}]
	\item[$i, \ j$] Index of buses.
	\item[$k$] Index of transmission elements.
	\item[$m$] Index of loads.
	\item[$n$] Index of generators.
	\item[$w$] Index of wind farms.
	\item[$t$] Index of load scenarios.
\end{IEEEdescription}

\subsection*{Variables}
\addcontentsline{toc}{section}{Nomenclature}
\begin{IEEEdescription}[\IEEEusemathlabelsep\IEEEsetlabelwidth{$V_1,V_2$}]
	\item[$P^g_{nt}$] Active power generation of generator $n$ in scenario $t$.
	\item[$P^g_{wt}$] Active power generation of wind farm $w$ in scenario $t$.
	\item[$P_{kt}$] Active power flow on branch $k$ in scenario $t$.
	\item[$\theta^p_{kt}$] Phase angle shift of the PST on branch $k$ in scenario $t$.
	\item[$\theta_{kt}$] The angle difference across the branch $k$ in scenario $t$.
	\item[$\delta_{k}$] Binary variable associated with placing a PST on branch $k$.
	\item[$\alpha_{k}$] Binary variable associated with line invesment for branch $k$.
\end{IEEEdescription} 

\subsection*{Parameters}
\addcontentsline{toc}{section}{Nomenclature}
\begin{IEEEdescription}[\IEEEusemathlabelsep\IEEEsetlabelwidth{$V_1,V_2$}]
	\item[$a_n^g$] Cost coefficient for generator $n$.
	\item[$P_{n}^{g,\min}$] Minimum active power output of generator $n$.
	\item[$P_{n}^{g,\max}$] Maximum active power output of generator $n$.
	\item[$P_{mt}^d$] Active power consumption of demand $m$ in scenario $t$.
	\item[$P_{wt}^a$] Available wind power of wind farm $w$ in scenario $t$.
	\item[$\theta_k^{p,\min}$] Minimum phase shift angle of PST on branch $k$.
	\item[$\theta_k^{p,\max}$] Maximum phase shift angle of PST on branch $k$.
    \item[$S_{kt}^{\max}$] Thermal limit of branch  $k$ in scenario $t$.
    \item[$A_k^p$] Investment cost of PST on branch $k$.
    \item[$\tilde{A}_k^p$] Annualized investment cost of PST on branch $k$.
    \item[$A_k^l$] Investment cost of line $k$.
    \item[$\tilde{A}_k^l$] Annualized investment cost of line $k$.
    \item[$A_p^{\max}$] Budget for investment on PST.
    \item[$A_l^{\max}$] Budget for investment one lines.
	\item[$N_{t}$] The number of operating hours for scenario $t$.
\end{IEEEdescription} 

\subsection*{Sets}
\addcontentsline{toc}{section}{Nomenclature}
\begin{IEEEdescription}[\IEEEusemathlabelsep\IEEEsetlabelwidth{$V_1,V_2$}]
	\item[$\mathcal{D}$] Set of loads.
	\item[$\mathcal{D}_i$] Set of loads located at bus $i$.
	\item[$\Omega_{L}$] Set of existing transmission lines.
	\item[$\Omega_{L}^{+}$] Set of prospective transmission lines.
	\item[$\Omega_{L}^i$] Set of transmission lines connected to bus $i$.
	\item[$\Omega_P$] Set of candidate lines to install PST.
	\item[$\Omega_{T}$] Set of scenarios.
	\item[$\mathcal{B}$] Set of buses.
	\item[$\mathcal{B}_{ref}$] Set of reference bus.
	\item[$\mathcal{G}$] Set of generators.
	\item[$\mathcal{G}_{i}$] Set of generators located at bus $i$.
	\item[$\mathcal{W}$] Set of wind farms.
	\item[$\mathcal{W}_{i}$] Set of wind farms located at bus $i$.
\end{IEEEdescription}

\section{Introduction}
\label{introduction}
\IEEEPARstart{T}{he} objective of the Transmission Expansion Planning (TEP) is to determine the best strategy to expand the existing power network in order to serve the growth of demand and generation in the future \cite{mybibb:tep_review}. The massive integration of wind power has introduced new challenges to the TEP problem. From the perspective of system planner, a rationally planned power network should not only improve the system reliability and security issues but also facilitate the integration of wind power \cite{mybibb:tep_es_conje}.   

The models and solution approaches for the TEP problem have been examined extensively in the technical literatures. Due to the nonlinear and non-convex characteristics of the power system planning problems, meta-heuristics methods \cite{mybibb:TEP_PSO,mybibb:TEP_DE,mybibb:yuan_chen_1} are usually adopted. These techniques have the advantage of model insensitivity and straightforward implementation. Nevertheless, the obtained solutions are not guaranteed to be global optimum and there is no indicator regarding the quality of the solution \cite{mybibb:TEP_huizhang2}. Reference \cite{mybibb:TEP_Heuristic_Compare} provides a review on the meta-heuristics methods in solving TEP problem.

Mathematical programming is another category of the solution techniques. The DC power flow model which considers only the active power and voltage angle is commonly leveraged \cite{mybibb:tep_large_wind,mybibb:TEP_letter}. The planning model is originally a mixed integer nonlinear program (MINLP). The product between a binary variable and a continuous variable can be linearized by introducing a disjunctive factor \cite{mybibb:TEP_2003}. Then the complete mode is transformed into a mixed integer linear program (MILP) which can be efficiently solved by commercial solvers. To incorporate the reactive power, several approximated AC TEP models have been proposed. The core idea of these models is to convexify the originally non-convex power flow equations. The techniques used for the convexification include: 1) linearizing the power flow equations around the operating points based on Taylor series \cite{mybibb:TEP_huizhang2,mybibb:TEP_Taylor}; 2) reformulating the power flow equations into semidefinite program (SDP) \cite{mybibb:TEP_Joshua} or second order cone program (SOCP) \cite{mybibb:Jabr_TEP}. The exact AC TEP model is difficult to be solved since there is no mature commercial solver to handle non-convex large scale MINLP. In \cite{mybibb:tep_exact_ac}, the AC TEP model is proposed and solved by combining interior point method and heuristic algorithm. To the best of our knowledge, the relaxed or exact AC TEP models are only utilized for small or medium scale system. 

To reduce the total planning cost, several technologies such as energy storage, transmission switching (TS) and series FACTS are introduced in the TEP process. In \cite{mybibb:tep_es_conje}, a long-term co-planning of transmission and energy storage considering wind power uncertainties is proposed. The transmission requirements and wind power curtailments can be reduced if the energy storages are appropriately placed in the system. In \cite{mybibb:tep_ts_denmark}, the authors evaluate the economic benefits of TS in the TEP problem. The problem is formulated as a two-stage stochastic model and solved by using branch-and-price algorithm. The authors in \cite{mybibb:tep_cvsr} present a security constrained multi-stage TEP model considering a continuously variable series reactor (CVSR). The CVSR serves as an additional corrective action during $N-1$ contingencies. In \cite{mybibb:tep_pst}, the phase shifting transformer (PST) is included in the TEP model and the genetic algorithm (GA) is utilized to solve the planning model. However, the wind power uncertainties are not considered.      
 
The investment decisions made by the system planner revolve around the electricity market. Thus, the market clearing conditions should be properly represented in the planning model. To achieve this goal, bilevel optimization model is commonly used in which the electricity market clearing is explicitly formulated in the lower level problem. A TEP model within a market environment is proposed in \cite{mybibb:bilevel_tep_conje}. In \cite{mybibb:tep_wind}, the authors jointly consider the transmission and wind farm expansion revolving around the electricity market by using bilevel model. In \cite{mybibb:tep_rpp}, with the linearized AC model, the authors propose a bilevel model to consider transmission and reactive power planning together. The load and wind uncertainties are represented by 5 scenarios in the single target planning year.  

This paper presents a bilevel optimization based TEP model under high penetration of wind power. The PST is introduced in the TEP process to provide extra flexibility. We adopt the static model which considers a single target year in the future \cite{mybibb:bilevel_tep_conje,mybibb:tep_wind}. The objective of the upper level problem is to minimize the total consumer payment, investment cost on new line and PST. The lower level problems represent a series of market clearing conditions under different load-wind scenarios. The contributions of this paper are twofold:
\begin{itemize}
	\item to propose a bilevel optimization based TEP model considering PST with significant wind power integration;
	\item to evaluate the benefits brought by the PST in the TEP process based on detailed case study results.
\end{itemize}   

The remaining sections are organized as follows. In Section \ref{reformulation}, the static model of PST is presented. Section \ref{optimization_model} provides details of the bilevel planning model. The solution approach to solve the bilevel optimization model is illustrated in Section \ref{solution_approach}. In Section \ref{case_studies}, the case studies based on IEEE 24-bus system is presented. Finally, some conclusions are given in Section \ref{conclusion}.

\section{Static Model of PST in DC Power Flow}
\label{reformulation}
\begin{figure}[!htb]
	\centering
	\includegraphics[width=0.3\textwidth]{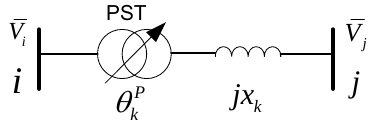}
	\caption{Static representation of PST in DCPF.}
	\label{pst_static}
\end{figure}

Fig. \ref{pst_static} depicts the steady state model of PST in DC power flow. It is modeled as a continuously variable phase angle $\theta_k^p$ in series with the transmission line reactance $x_k$.     

The power flow on the transmission line can be expressed as:
\begin{align}
&P_k=b_k(\theta_k+\delta_k\theta_k^p) \label{pf_pst} \\
&\theta_k^{p,\min} \le \theta_k^p \le \theta_k^{p,\max}  \label{angle_range}
\end{align}
where 
\begin{equation}
b_{k}=\frac{1}{x_k}   \\
\end{equation}

In (\ref{pf_pst}), a binary variable $\delta_k$ is used to flag the installation of the PST on transmission line $k$. It can be seen that constraint (\ref{pf_pst}) is nonlinear due to the product between $\delta_k$ and $\theta_k^p$. To linearize the nonlinear term, we introduce an auxiliary variable $\psi_k^p$ and constraints (\ref{pf_pst}) and (\ref{angle_range}) can be reformulated as \cite{mybibb:tcpst_milp}:
\begin{align}
&P_k=b_k\theta_k+b_k\psi_k^p \label{pst_refor}  \\
&\delta_k\theta_k^{p,\min} \le \psi_k^p \le \delta_k\theta_k^{p,\max} 
\end{align}

If transmission line $k$ is not selected to install a PST, i.e., $\delta_k=0$, $\psi_k^p$ will be zero. Otherwise, $\psi_k^p$ will be bounded by the maximum and minimum phase shift angle. 

\section{Problem Formulation}
\label{optimization_model}
As mentioned in the introduction section, the planning problem can be formulated as a bilevel optimization model. The objective of the upper level problem is to minimize the total consumer payment in the target planning year and the investment cost on new line and PST. The market clearing conditions under different load-wind scenarios are considered in the lower level problems. The complete optimization model is given as:
\begin{subequations}
	\label{bilevel_entire}
	\begin{align}
	&\min_{\Xi_{\text{UL}}} \ \ \sum_{k\in \Omega_P}\tilde{A}_k^p\delta_k+\sum_{k\in \Omega_L^+}\tilde{A}_k^l\alpha_k \nonumber \\
	&\ \ \ \ \ \ \ \ \ +\sum_{t\in \Omega_T}N_t \sum_{i\in \mathcal{B}}\lambda_{it}\sum_{m \in \mathcal{D}_i}P_{mt}^d \label{objective_upper}  \\
	&\text{subject to}  \nonumber \\
	&\sum_{k \in \Omega_P}A^p_k \delta_k \le A_p^{\max}  \label{budget_p}  \\
	&\sum_{k \in \Omega_L^{+}}A^l_k \alpha_k \le A_l^{\max}  \label{budget_line}  \\
	&\text{where} \  \lambda_{it} \in \arg \{  \nonumber \\
	&\min_{\Xi_{\text{LL}}} \ \ \sum_{n\in \mathcal{G}}a_{n}^gP^g_{nt} \label{obj_lower} \\
	&\text{subject to}  \nonumber \\
	&\sum_{n \in \mathcal{G}_i}P^g_{nt}+\sum_{w \in \mathcal{W}_i}P^g_{wt}-\sum_{m \in \mathcal{D}_i}P^d_{mt} \nonumber  \\
	&=\sum_{k\in \Omega_{L}^i } P_{kt}:  \lambda_{it}, \ i\in \mathcal{B}  \label{p_bal}  \\
	&P_{kt}=b_k \theta_{kt}:\phi^l_{kt}, \ k\in \Omega_{L}\backslash \Omega_{P} \label{norm_line} \\
	&P_{kt}=b_k(\theta_{kt}+\psi_{kt}^P):\phi^p_{kt}, \ k\in \Omega_{P} \label{pst_line} \\
	&\delta_k\theta_k^{p,\min} \le \psi_{kt}^p \le \delta_k\theta_k^{p,\max}:\phi_{kt}^{p^{\min}},\phi_{kt}^{p^{\max}}, \ k\in \Omega_P \label{pst_anlge_lim}   \\
	&P_{kt}-b_k\theta_{kt}+M_k(1-\alpha_k) \ge 0:\phi_{kt}^{+^{\min}}, \ k \in \Omega_L^+  \label{prospective_1}  \\
	&P_{kt}-b_k\theta_{kt}-M_k(1-\alpha_k) \le 0:\phi_{kt}^{+^{\max}}, \ k \in \Omega_L^+  \label{prospective_2}  \\
 	&P_{n}^{g,\min}\le P_{nt}^g \le P_{n}^{g,\max}:\xi_{nt}^{g^{\min}},\xi_{nt}^{g^{\max}} \label{Pg_limit} \\ 
	&0 \le P_{wt}^g \le P_{wt}^a:\xi_{wt}^{g^{\min}},\xi_{wt}^{g^{\max}}  \label{wind_limit}    \\
	& -S_{kt}^{\max} \le P_{kt} \le S_{kt}^{\max}:\eta_{kt}^{l^{\min}},\eta_{kt}^{l^{\max}}, \ k\in \Omega_L \label{Slimit_E} \\
	&-\alpha_kS_{kt}^{\max} \le P_{kt} \le \alpha_kS_{kt}^{\max}:\eta_{kt}^{+^{\min}},\eta_{kt}^{+^{\max}}, \ k\in \Omega_L^+ \label{Slimit_C} \\
	& -\pi \le \theta_{it} \le \pi:\gamma_{it}^{\min},\gamma_{it}^{\max}, \ i\in \mathcal{B}/\mathcal{B}_{ref} \label{bus_angle_lim}   \\
	& \theta_{it}=0: \gamma_{it}^{ref}, \ i \in \mathcal{B}_{ref}     \label{ref_angle} \ \  \}
	\end{align}
\end{subequations}   
Constraints (\ref{p_bal})-(\ref{ref_angle}) hold $\forall n \in \mathcal{G}, m \in \mathcal{D}, w\in \mathcal{W}, t \in \Omega_T$. The dual variables associated with the constraints in the lower level problem are provided following a colon. The optimization variables in the lower level problem are: $\splitatcommas{\Xi_{\text{LL}}=\{P^g_{nt},P^g_{wt},P_{kt},\psi_{kt}^P,\theta_{kt},\theta_{it},\lambda_{it},\phi_{kt}^l,\phi^p_{kt},\phi_{kt}^{p^{\min}},\phi_{kt}^{p^{\max}},\phi_{kt}^{+^{\min}}, \phi_{kt}^{+^{\max}},\xi_{nt}^{g^{\min}},\xi_{nt}^{g^{\max}},\xi_{wt}^{g^{\min}},\xi_{wt}^{g^{\max}},\eta_{kt}^{l^{\min}},\eta_{kt}^{l^{\max}},\gamma_{it}^{\min},\gamma_{it}^{\max},\gamma_{it}^{ref} \}}$. The optimization variables in the upper level problem are: $\Xi_{\text{UL}}=\{\delta_k,\alpha_k,\Xi_{\text{LL}}\}$.

The objective function of the upper level problem, i.e., (\ref{objective_upper}), comprises three terms. Specifically, the first two terms denote the annualized investment cost on PST and transmission line. The third term is the total consumer payment. It is assumed that the price of consumer payment is the locational marginal price (LMP) at the bus where the load is located \cite{mybibb:tep_wind}. The LMP for scenario $t$ at bus $i$ is the dual variable associated with the power balance constraint (\ref{p_bal}), i.e., $\lambda_{it}$. Then the consumer payment for each scenario $t$ at bus $i$ is the corresponding demand multiplied by $\lambda_{it}$. Constraints (\ref{budget_p}) and (\ref{budget_line}) impose budgets on the investment of PST and transmission line respectively. 

The lower level problems designate market clearing conditions under different load-wind scenarios. The objective of each lower level problem is to maximize the social welfare. Consider a perfect inelastic load, maximizing the social welfare is equivalent to minimizing the production cost, which is denoted by (\ref{obj_lower}). We assume that the marginal cost for wind power is zero. The active power balance is enforced by constraint (\ref{p_bal}). Constraint (\ref{norm_line}) defines the power flow through the existing lines without PST. Constraints (\ref{pst_line}) and (\ref{pst_anlge_lim}) indicate the power flow on the candidate lines to install PST. The active power flow through the prospective lines are provided in constraints (\ref{prospective_1}) and (\ref{prospective_2}). When the prospective line is selected to be built, i.e., $\alpha_k=1$, the line flow equations are enforced to hold; when the transmission line is not selected, the sufficiently large number $M_k$ will ensure that the two constraints are not binding and become redundant \cite{mybibb:tep_cvsr}. Constraint (\ref{Pg_limit}) enforces the generation limits for the conventional generator. The dispatched wind power is limited by its available amount in (\ref{wind_limit}). This constraint also allows wind power curtailment if needed. Constraints (\ref{Slimit_E}) and (\ref{Slimit_C}) represent the thermal limits for the existing and prospective transmission lines. The bus angle bounds are provided in constraint (\ref{bus_angle_lim}). Finally, the bus angle of the reference bus is set to be zero in constraint (\ref{ref_angle}).

\section{Solution Approach}
\label{solution_approach}
With $\delta_k$ and $\alpha_k$ from the upper level problem, the lower level problems are pure linear program (LP). Therefore, the bilevel model can be transformed into a mathematic program with equilibrium constraints (MPEC) model by replace each lower level problem with its KKT conditions or primal-dual formulation \cite{mybibb:tep_wind,mybibb:xinfang_bilevel2}. The KKT based formulation includes a large number of complementary constraints which are linearized by introducing auxiliary binary variables \cite{mybibb:tep_wind}. Thus, we leverage the primal-dual based formulation which is computationally more friendly than the KKT based one \cite{mybibb:mip_bilevel}.   

We define $\bm{x}$ as the upper level binary variables and $\bm{y_t}$ as the lower level primal variables for scenario $t$. The bilevel optimization model (\ref{bilevel_entire}) can be compactly written as:
\begin{subequations}
	\label{bilevel_entire_comp}
	\begin{align}
	&\min_{\Xi_{\text{UL}}} \ \ \bm{f}^T\bm{x}+\sum_{t \in \Omega_T}N_t\bm{g_t}^T\bm{\lambda_t}  \label{comp_obj} \\
	&\bm{Ax} \le \bm{b}    \label{com_upper_c}   \\
	&\text{where} \ \bm{\lambda_t} \in \arg \{ \ \ \min_{\Xi_{\text{LL}}} \ \ \bm{w_t}^T\bm{y_t}  \label{com_low_obj} \\
	&\text{s.t.} \ \bm{Py_t}\le \bm{r_t}-\bm{Kx}:\bm{\mu_t} \label{com_low_ineq}  \\
	&\ \ \ \ \bm{Ey_t}=\bm{h_t}:\bm{\lambda_t}  \ \  \} \ \ \forall t \in \Omega_T \label{com_low_eq}
	\end{align} 
\end{subequations}
In (\ref{bilevel_entire_comp}), $\bm{\lambda_t}$ and $\bm{\mu_t}$ comprise all the dual variables associated with the equality and inequality constraints in scenario $t$.

In primal-dual reformation, each of the lower level problem should be replaced by its primal constraints, dual constraints and the strong duality. The bilevel model (\ref{bilevel_entire_comp}) can be transformed into a single level model as follows:
\begin{subequations}
	\label{bilevel_single}
	\begin{align}
	&\min_{\Xi_{\text{UL}}} \ \ \bm{f}^T\bm{x}+\sum_{t \in \Omega_T}N_t\bm{g_t}^T\bm{\lambda_t}  \label{comp_obj_s} \\
	&\bm{Ax} \le \bm{b}    \label{com_upper_c_s}   \\
	&\bm{Py_t}\le \bm{r_t}-\bm{Kx} \label{com_low_ineq_s}  \\
	&\bm{Ey_t}=\bm{h_t} \label{com_low_eq_s}  \\
	&\bm{P}^T\bm{\mu_t}+\bm{E}^T\bm{\lambda_t}+\bm{w_t}=0  \label{dual_eq}  \\
	&\bm{\mu_t} \ge 0 \label{dual_ineq}  \\
	&\bm{w_t}^T\bm{y_t}=(\bm{x}^T\bm{K}^T-\bm{r_t}^T)\bm{\mu_t}-\bm{h_t}^T\bm{\lambda_t} \label{strong_dual} 
	\end{align} 
\end{subequations}
 
Constraint (\ref{com_upper_c_s})-(\ref{strong_dual}) hold $\forall t \in \Omega_T$. Constraints (\ref{com_low_ineq_s}) and (\ref{com_low_eq_s}) denote the primal constraints. The dual constraints are provided in (\ref{dual_eq}) and (\ref{dual_ineq}). Constraint (\ref{strong_dual}) represents the strong duality theorem, i.e., the objective of the primal and dual problem should be equal if the optimization model is convex. 

Note that in constraint (\ref{strong_dual}), there is a bilinear term which is the product between the binary variable $\bm{x}$ and continuous variable $\bm{\mu_t}$. This bilinear term can be easily linearized using the big-M method \cite{mybibb:svc_gm_2017,mybibb:Tao_bigM}. Thus, the optimization model (\ref{bilevel_single}) is an MILP which can be efficiently solved by commercial solvers.  

\section{Numerical Case Studies}
\label{case_studies}
We test our proposed planning model on the IEEE 24-bus system. We leverage the MATLAB based toolbox YALMIP \cite{mybibb:YALMIP} to model the complete MILP problem and the CPLEX solver \cite{mybibb:CPLEX} to solve the model. The computer used to perform all the simulations has an Inter(R) Core(TM) i5-6300U CPU
@ 2.40 GHz and 8.00 GB of RAM. 

The investment cost on the transmission line can be estimated by its length, capital cost per mile and the cost multiplier \cite{mybibb:wecc_cost}. For the PST, the investment cost is based on the current rating of the corresponding transmission line \cite{mybibb:GA_FACT_market}. The cost coefficient of PST is selected to be 100 \$/kVA and we allow the phase shift angle to vary from $-5^{\circ}$ to $+5^{\circ}$ \cite{mybibb:pst_cost1}. The annualized investment costs on transmission line and PST are converted from their total costs using the following equation \cite{mybibb:PSO_SQP_FACTS,mybibb:investment_naps}:
\begin{equation}
\tilde{A}_k={A}_k \cdot \frac{d(1+d)^{LT}}{(1+d)^{LT}-1}  \label{cost_convert}
\end{equation}    
In (\ref{cost_convert}), $LT$ is the life time of the new facilities and $d$ is the interest rate. In this work, the life time of the transmission line and PST is assumed to be 20 and 15 years respectively. The interest rate $d$ is chosen to be 5\%.

The sufficiently large constant $M_k$ in constraint (\ref{prospective_1}) and (\ref{prospective_2}) is selected to be $|2\pi b_k|$. Finally, the ``mipgap" in CPLEX is set to be 0.1\%. 

\subsection{Load and Wind Data}
In the simulations, the normalized electric load profile from the 2015 ISO New England hourly demand reports \cite{mybibb:load_data} is used to represent the annual load for the test system. Moreover, the hourly wind power capacity factor for the target year is assumed to follow the profile given by \cite{mybibb:wind_data} with the wind turbine model to be GE 1.5sl \cite{mybibb:ge_wind_turbine} and location to be the North part of Denmark. We then use \textsl{K-means} method to reduce the the number of scenarios from 8760 to 10. The reason we leverage \textsl{K-means} for the scenario reduction is that it allows to retain the correlations between the demand and wind power \cite{mybibb:tep_wind}. The number of operating hours, load levels and wind capacity factors for each of the 10 scenarios are provided in Table \ref{load_wind_s}.

\begin{table}[!htb]
	\centering
	\caption{Load and Wind Scenarios}
	\label{load_wind_s}
	\begin{tabular}{c c c c c c }
		\hline
		Scenario \#&1&2&3&4&5  \\
		\hline
		Load level&0.8307&0.5456&0.5220&0.6999&0.7301   \\
		\hline
		Wind capacity factor&0.4287&0.7280&0.0946&0.7739&0.1523  \\
		\hline
		Number of hours&355&742&1323&553&927   \\
		\hline
		Scenario \#&6&7&8&9&10  \\
		\hline
		Load level&0.5224&0.6496&0.4999&0.5556&0.6713   \\
		\hline
		Wind capacity factor&0.5454&0.3616&0.3577&0.2185&0.5659  \\
		\hline
		Number of hours&780&1057&900&1328&795   \\
		\hline
	\end{tabular}
\end{table}

\subsection{IEEE 24-Bus System}
The IEEE 24-bus system has 24 buses, 38 transmission lines and 32 generators. The system data can be found in \cite{mybibb:MATPOW}. To create congestions, the peak loads and generator capacities are 1.5 times the values given in \cite{mybibb:MATPOW}. Moreover, the thermal limits of the transmission lines are decreased to 60\% of the original values. Two wind farms with the capacity of 1200 MW are located at bus 10 and 14. It is assumed that the wind capacity factors of wind farm at bus 14 follow the values provided in Table \ref{load_wind_s} and the capacity factors of wind farm at bus 10 are 10\% lower than those values. 

We consider 7 transmission corridor as the candidate locations to install the transmission lines. The data for the prospective lines is provided in Table \ref{line_data}. The sensitivity approach in \cite{mybibb:PSO_SQP_FACTS} is leveraged to select 10 existing lines to install PST. 

\begin{table}[!htb]
	\centering
	\caption{Prospective Line Data}
	\label{line_data}
	\begin{tabular}{c c c c c}
		\hline
		From&To&Reactance&Capacity&Investment cost   \\
		bus&bus&(p.u.)&(MW)&(M\$)  \\
		\hline
		1&2&0.0139&105&3.9094  \\
		\hline
		2&6&0.1920&105&54.0000  \\
		\hline
		6&10&0.0605&105&17.0156  \\
		\hline
		7&8&0.0614&105&17.2688  \\
		\hline
		8&9&0.1651&105&46.4344  \\
		\hline
		8&10&0.1651&105&46.4344  \\
		\hline
		9&12&0.0839&240&132.4737  \\
		\hline
	\end{tabular}
\end{table}
	
\begin{table*}[t]
	\centering
	\caption{IEEE 24-Bus System Results for Different Cases}
	\label{ieee_24_results}
	\begin{tabular}{c c c c c c c c c c c}
		\hline
		&\multicolumn{2}{c}{Wind curtailment }&\multirow{3}{1.2 cm}{\centering Lines built}&\multirow{2}{1.0 cm}{\centering Investment on lines}&\multirow{3}{1.2 cm}{\centering PST locations}&\multirow{2}{1.0 cm}{\centering Investment on PSTs}&\multirow{2}{1.0 cm}{\centering Consumer payment}&\multirow{2}{1.0 cm}{\centering Objective value}&\multirow{2}{1.2 cm}{\centering Wind penetration}&\multirow{2}{1.0 cm}{\centering Time (s)}   \\
		&\multicolumn{2}{c}{($\times 10^6$ MWh)}&&on lines&&&&&&   \\
		\cline{2-3} 
		&10&14&&(M \$)&&(M \$)&(M \$)&(M \$)&level (\%)&   \\
		\hline
		C1&0.5930&0.3335&-&-&-&-&430.3055&430.3055&28.8426&0.2471   \\
		\hline
		&&&6-10,7-8&&&&&&&   \\
		C2&0.0102&0.4920&8-9,8-10&20.8331&-&-&348.3527&369.1858&30.7260&1.3148   \\
		&&&9-12&&&&&&&    \\
		\hline
		\multirow{2}{*}{C3}&\multirow{2}{*}{0}&\multirow{2}{*}{0.4057}&6-10,8-9&\multirow{2}{1.0 cm}{\centering 19.4474}&\multirow{2}{1.2 cm}{\centering 3-9}&\multirow{2}{1.0 cm}{\centering 1.0116}&\multirow{2}{1.0 cm}{\centering 295.7458}&\multirow{2}{1.0 cm}{\centering 316.2048}&\multirow{2}{1.2 cm}{\centering 31.1456}&\multirow{2}{1.0 cm}{\centering 2.1668}   \\
		&&&8-10,9-12&&&&&&&   \\
		\hline
		\multirow{2}{*}{C4}&\multirow{2}{*}{0.0872}&\multirow{2}{*}{0.2506}&1-2,2-6&\multirow{2}{1.0 cm}{\centering 12.0988}&\multirow{2}{1.2 cm}{\centering 1-5,3-9}&\multirow{2}{1.0 cm}{\centering 2.0232}&\multirow{2}{1.0 cm}{\centering 295.4857}&\multirow{2}{1.0 cm}{\centering 309.6077}&\multirow{2}{1.2 cm}{\centering 31.4556}&\multirow{2}{1.0 cm}{\centering 3.3508}   \\
		&&&8-9,8-10&&&&&&&   \\
		\hline
		
	\end{tabular}
\end{table*}    

Four cases are considered in the simulations:
\begin{itemize}
	\item Case 1: The budgets for the transmission lines and PST are both zero, i.e., no reinforcement is considered.
	\item Case 2: The budget for the transmission lines is unlimited while the budget of PST is zero.
	\item Case 3: The budget for the transmission lines is unlimited while the budget of PST is \$15 million.
	\item Case 4: The budget for the transmission lines is unlimited while the budget of PST is \$30 million.
\end{itemize}    

The simulation results for different cases are provided in Table \ref{ieee_24_results}. The second and third columns provide the annual wind power curtailment for each wind farm. The fourth and fifth columns give the lines to be built and the annualized investment cost on the lines. The sixth and seventh columns indicate the lines to be installed with the PST and annualized investment cost on PST. The consumer payment, i.e., the third term in (\ref{objective_upper}), is given in the eighth column. The value of the objective function is provided in the ninth column. The tenth column indicates the wind penetration level which is defined as the percentage of load that can be covered by wind power on an annual basis. Finally, the computational time is provided in the last column.  

From Table \ref{ieee_24_results}, it can be seen that the consumer payment without any network reinforcement is \$430.31 M and the wind penetration level is 28.84\%. When only TEP is allowed, 5 transmission lines with the annualized investment cost of \$20.83 M are built. This reinforcement decreases the total consumer payment to \$348.35 M. Moreover, the wind penetration level is increased by 1.88\% compared to the value in Case 1. The introduction of PST brings the economic benefits to the TEP process. In Case 3, 1 PST is selected to be installed on line 3-9 and only 4 transmission lines are built. The construction of line 7-8 is avoided. Although the difference between the investment of Case 2 and Case 3 is trivial, a significant reduction on the consumer payment is observed. In addition, the wind penetration level is increased to 31.15\%. When the budget for the PST is \$30 M, i.e., Case 4, the planning model suggests to build 4 transmission lines and 2 PSTs. The construction of the most expensive line 9-12 is avoided, which reduces the investment cost to \$14.12 M. When compared the value of objective function in Case 2 and Case 4, it can be observed that a total savings of \$60 M is achieved with the installation of PSTs. 

\begin{figure}[!htb]
	\centering
	\includegraphics[width=0.45\textwidth]{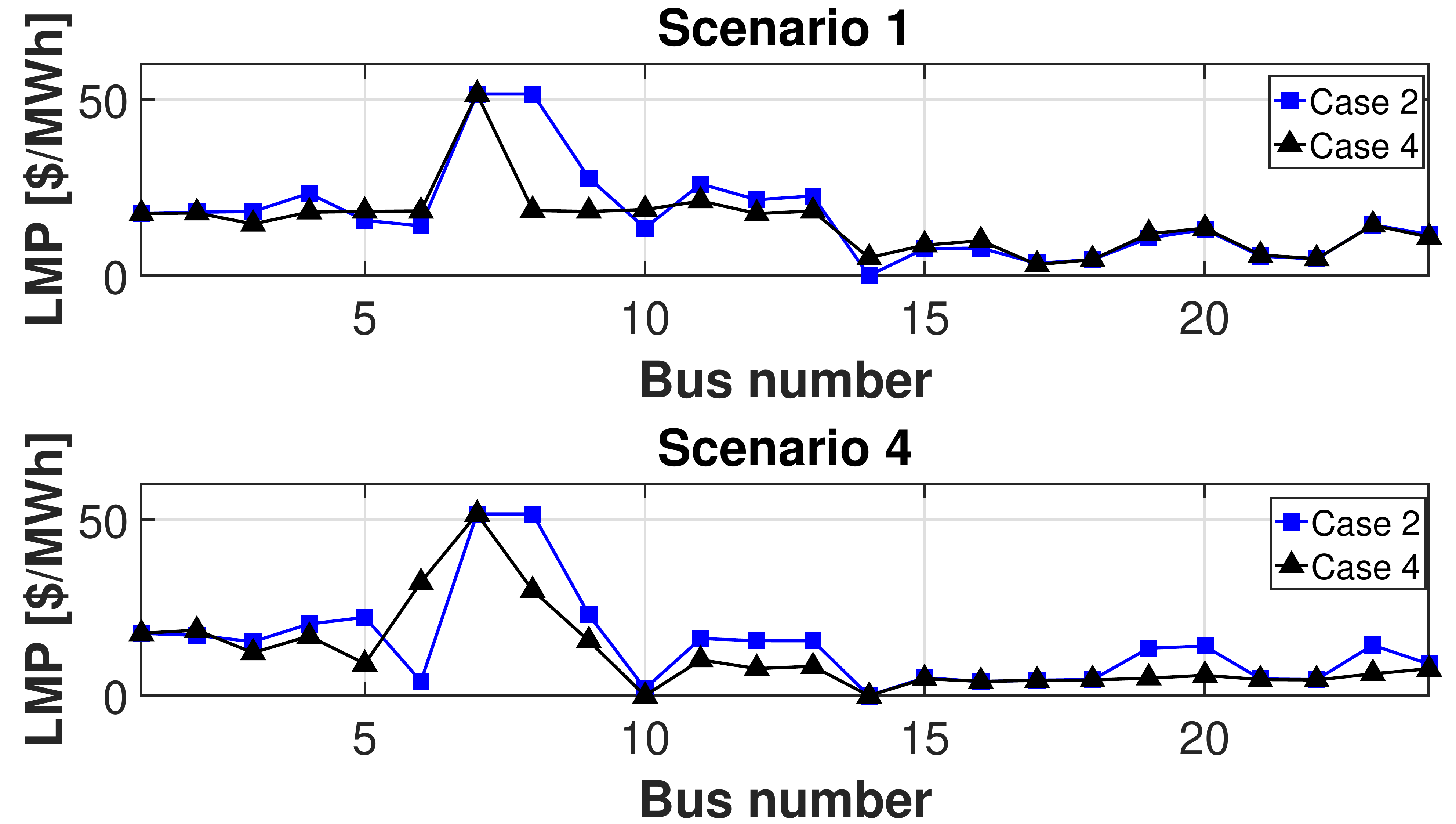}
	\caption{Comparison of the LMPs at each bus.}
	\label{LMP}
\end{figure}

Fig. \ref{LMP} depicts the LMPs for 2 scenarios in Case 2 and Case 4. From Table \ref{load_wind_s}, it can be seen that scenario 1 represents the scenario with the highest load level and scenario 2 indicates the scenario with the highest wind level. As can be observed from Fig. \ref{LMP}, at the majority buses, the LMP in Case 4 is lower than that in Case 2 for both scenarios. The installation of PSTs tends to relieve the congestion and results in more flat LMPs.

\begin{figure}[!htb]
	\centering
	\includegraphics[width=0.3\textwidth]{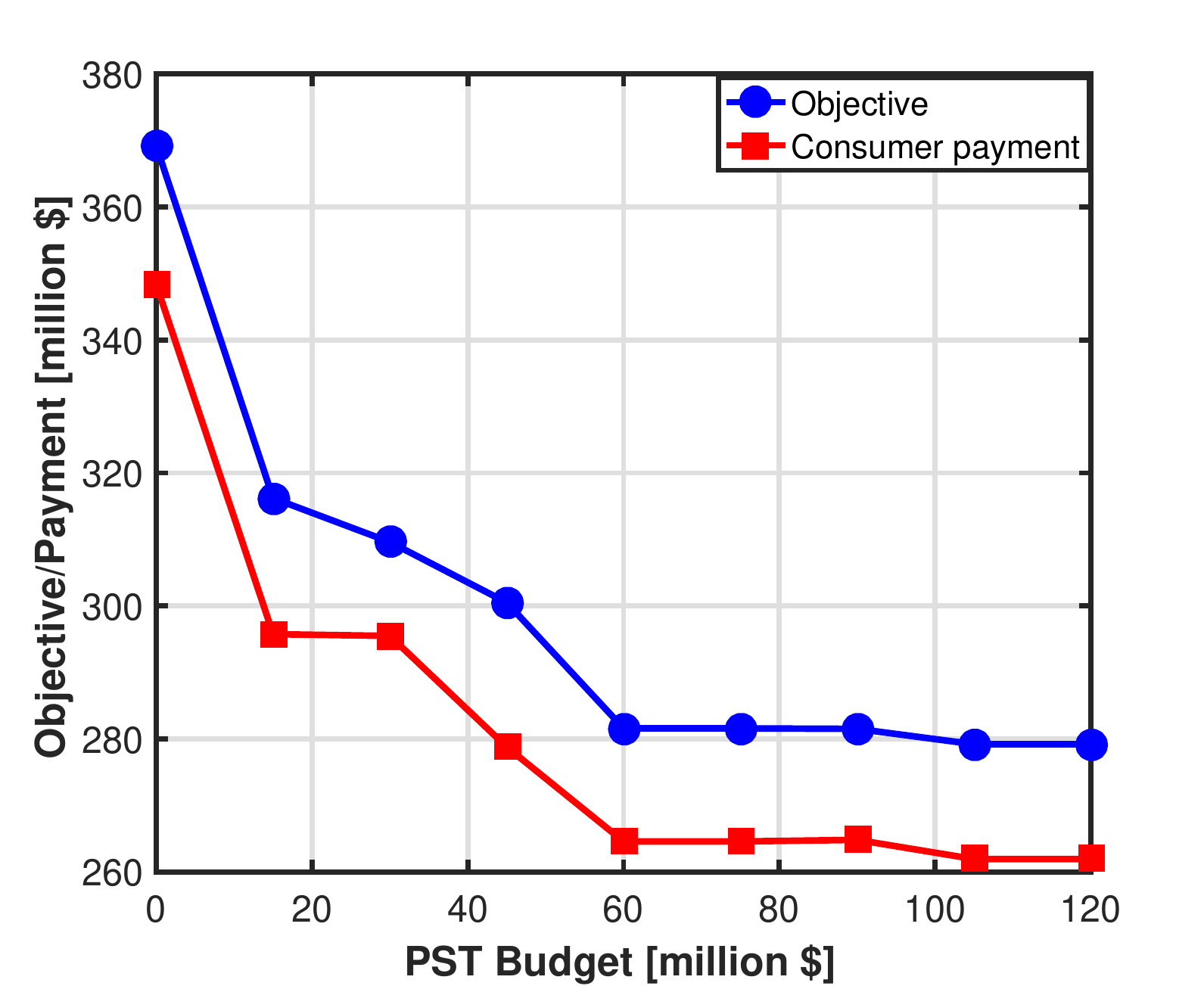}
	\caption{Optimal value of the consumer payment and the objective function for different budgets on PST.}
	\label{iter_obj}
\end{figure}

Fig. \ref{iter_obj} illustrates the objective value and consumer payment as a function of the budget on PST. Note that the budget for the transmission lines is assumed unlimited for this plot. As can be seen from the figure, both of the consumer payment and the objective value decreases as the budget on PST increases. However, most of the cost reductions are achieved when the budget on PST is \$30 M or \$45 M, i.e., 2 or 3 PSTs. The objective value can be future reduced if more PSTs are introduced but the amount of the cost reduction is small.     

\section{Conclusions}
\label{conclusion}
 This paper proposes a bilevel optimization based TEP model considering PST under high penetration of wind power. The proposed model seeks to identify the investment strategy on transmission line and PST within a market environment. The load-wind uncertainties are represented by a number of scenarios. To transform the bilevel model into a single level problem, each lower level problem is replaced by its corresponding primal-dual formulation. The numerical results on IEEE 24-bus system illustrate that the introduction of PST in the TEP problem not only allows reduced total planning cost, but also facilitates the wind power integration.

\bibliographystyle{IEEEtran}
\bibliography{IEEEabrv,mybibb}

\end{document}